\documentclass[12pt]{article}
\usepackage{amssymb}
\usepackage{amsmath}
\newtheorem{teo}{Theorem}

\newtheorem{cor}{Corollary}

\title{Matched wavelets for equidistant points
}
\author{
E. A. Lebedeva
\footnote{Mathematics and Mechanics Faculty, Saint Petersburg State University,
Universitetsky prospekt, 28, Peterhof,  Saint Petersburg, 198504, Russia
 }
}
\date{
 ealebedeva2004@gmail.com}

\begin{document}
\maketitle

\newcommand{\nul}{{\bf0}}
\newcommand{\rd}{{\mathbb R}^d}
\newcommand{\zd}{{\mathbb Z}^{d}}
\renewcommand{\r}{{\mathbb R}}
\newcommand{\z} {{\mathbb Z}}
\newcommand{\cn} {{\mathbb C}}
\newcommand{\n} {{\mathbb N}}

\begin{abstract}
Matched wavelets interpolating equidistant  data  are designed. These wavelets form   Riesz bases. Meyer wavelets that interpolate data on a particular uniform lattice are found.       
\end{abstract}

\textbf{Keywords} matched wavelet, Meyer wavelet, interpolation wavelet, Riesz basis. 

\textbf{AMS Subject Classification}:  42C40

\section{Introduction}	

In many topics  of signal processing, wavelets that are matched to a signal of interest are required. To this end different optimization techniques are often employed. For example, in \cite{ChapaRao}, the authors find a bandlimited wavelet that minimizes $L_2$-distance between a wavelet and a desired signal in the Fourier domain. In \cite{KarelPeeters}, to design orthogonal wavelets  with compact support and vanishing moments the following optimization criteria are used: sparsity by $L_1$-minimization and by $L_4$-norm maximization.  

In this paper, we present another approach to the question.     
Under a matched wavelet we understand a wavelet function interpolating data on some lattice. More precisely, given two sequences $\alpha_k$, $\gamma_k \in \mathbb{R}$, $k\in\mathbb{Z}$. A wavelet function $\psi$ is called \textit{a matched wavelet} on the lattice $\alpha_k$, $k \in \mathbb{Z},$  if $\psi(\alpha_k) = \gamma_k$ for any $k\in \mathbb{Z}.$

The paper is organized as follows. In section \ref{Not}, we recall necessary definitions and notations. In section \ref{main}, we  design a matched  wavelet function that generates  Riesz basis and interpolates very wide class of data on the lattice $\alpha_n = n+1/2,$ $n\in\mathbb{Z}$ (Theorem \ref{riesz}). Additionally, the case of finite data is discussed. For a particular choice of the lattice (namely, $\alpha_n = 3n+1/2,$ $n\in\mathbb{Z}$) we find an orthogonal  wavelet function (the Meyer wavelet) that interpolates wide class of data on the lattice (Theorem \ref{Mey}).

\section{Notations and Auxiliary Results}
\label{Not}
 Recall the definition of Meyer wavelets \cite{NPS}. By $\theta(\xi)$ denote an odd continuously differentiable  function equal to $\pi/4$ for $\xi>\pi/3.$ By $\lambda(\xi)$ denote an even function defined on $[0,\,\infty)$ as follows
$$
\lambda(\xi) = \left\{
\begin{array}{ll}
\pi/4 + \theta(\xi-\pi), & \xi \in [2\pi /3,\, 4\pi/3], \\
\pi/4 - \theta(\xi/2-\pi), & \xi \in [4\pi /3,\, 8\pi/3], \\
0, & \xi \in [0,\, 2 \pi/3)\cup (8 \pi/3,\,\infty). 
\end{array}
\right.
$$
The Fourier transform of Meyer scaling function is equal to
$$
\widehat{\varphi^M}(\xi) = \left\{
\begin{array}{ll}
\cos \lambda(\xi), & |\xi|\leq 4 \pi /3, \\
0,& |\xi| > 4 \pi /3,
\end{array}
\right.
$$
where we choose the Fourier transform in the form $\displaystyle \widehat{f}(\xi)=\int_{\mathbb{R}}f(t){\rm e}^{-{\rm i} \xi t}\, dt.$ 
Corresponding mask $m^M(\xi)$ is equal to $\widehat{\varphi^M}(2\xi)$ as $|\xi|<\pi$.   
Then the Fourier transform of Meyer wavelet is defined by
$$
\widehat{\psi^M}(\xi) = {\rm e}^{-{\rm i}\xi/2} \overline{m^M(\xi/2+\pi)} \widehat{\varphi^M}(\xi/2) = {\rm e}^{-{\rm i}\xi/2} \sin \lambda (\xi)$$
or in the time domain 
$$
\psi^M(t)=\frac{1}{2\pi} \int_{\mathbb{R}} \cos((t-1/2)\xi) \sin\lambda(\xi)\,d\xi.
$$
We explore the traditional notation for the dilations and shifts of a function 
$f_{j,k}(x):=2^{j/2}f(2^j x+k),$ $j,k\in\mathbb{Z}.$

\section{Results}
\label{main}

In Theorem \ref{riesz} we design a wavelet function that generates  Riesz basis and interpolates very wide class of data on the lattice $n+1/2,$ $n\in\mathbb{Z}.$

\begin{teo}
\label{riesz}
Suppose $\gamma_k \in \mathbb{C},$ $k \in \mathbb{Z},$ is a sequence such that 
\begin{equation}
\label{Gam}
\exists A,B >0 \ \  A\leq \left|\Gamma(\xi)\right|^2\leq B,
\end{equation}
where $\Gamma(\xi) = \sum_{k\in \mathbb{Z}} \gamma_k {\rm e}^{{\rm i}k \xi}$. 
Then there exists a wavelet Riesz basis $\psi_{j,k},$ $j,k\in\mathbb{Z}$ of $L_2(\mathbb{R})$  with bounds $A$ and $B$ such that 
$$
\psi(k+1/2)=\gamma_k, \ \  k\in\mathbb{Z}.
$$ 
\end{teo}

\textbf{Proof. } Let $\varphi^I$ be a scaling function of an orthogonal multiresolution analysis,  and $\varphi^I(k)=0$ for any $k\in\mathbb{Z}\setminus\{0\}$, $\varphi^I(0)=1$. In the literature, such a function is sometimes called  interpolating or cardinal \cite{AlUn}.   Suppose, that a corresponding wavelet function $\psi^I$ generates an orthonormal basis $\psi^I_{j,k},$ $j,k\in\mathbb{Z}$, of $L_2(\mathbb{R}).$ By $W_j(\psi^I):=\overline{{\rm span}\{\psi^I_{j,k}\ :\ k\in\mathbb{Z}\}},$ $j \in \mathbb{Z}$, we denote wavelet spaces generated by $\psi^I.$ 
Let us  define 
\begin{equation}
\label{psi1}
\psi(x) = \sum_{k\in\mathbb{Z}} \gamma_k \psi^{I}(x+k).
\end{equation}
The function $\psi^I$ inherits the interpolation property from the scaling function $\varphi^I$, namely, $\psi^I(n+1/2)=0,$ if $n\neq 0,$ and  $\psi^I(1/2)=1,$ so 
$$
\psi\left(n+\frac{1}{2}\right) = \sum_{k\in\mathbb{Z}} \gamma_k \psi^I\left(n+ k +\frac{1}{2}\right) = \gamma_n. 
$$ 

Since (\ref{psi1}) can be  rewritten in the Fourier domain as $\widehat{\psi}(\xi) = \Gamma(\xi)\widehat{\psi^I}(\xi)$ and $\displaystyle \left|\Gamma(\xi)\right| \leq \sqrt{B},$ it follows that $\Gamma \in L_2[-\pi,\,\pi]$, therefore, $(\gamma_k)_k \in l_2,$ thus $\psi \in W_0(\psi^I)$.

Using orthogonality of the shifts $\psi^I(x+k),$ $k\in\mathbb{Z}$, we note that 
$$
\sum_{k\in \mathbb{Z}} \left|\widehat{\psi}(\xi + k)\right|^2 = \left|\Gamma(\xi)\right|^2 
\sum_{k\in \mathbb{Z}} \left|\widehat{\psi^I}(\xi + k)\right|^2 = \left|\Gamma(\xi)\right|^2. 
$$     
So, it follows from (\ref{Gam}), that $\psi(x+k),$ $k\in\mathbb{Z},$ is a Riesz system. Therefore, \cite[Theorem 1.1.2]{NPS}  $\psi(x+k),$ $k\in\mathbb{Z},$ is a Riesz basis in 
$$
W_0(\psi):=\left\{f(\xi) = \sum_{k\in\mathbb{Z}}c_k\psi(\xi+k) \ : \ (c_k)_k \in l_2 \right\}
=\overline{{\rm span}\left\{\psi_{0,k}\ : \ k\in\mathbb{Z}\right\}}.
$$
Denote  $W_j(\psi)= \overline{{\rm span}\left\{\psi_{j,k}\ : \ k\in\mathbb{Z}\right\}},$
 $j\in\mathbb{Z}.$

Since $\widehat{\psi^I}(\xi) = \left(\Gamma(\xi)\right)^{-1}\widehat{\psi}(\xi)$ and $\displaystyle \left|\Gamma(\xi)\right|^{-2} \leq A^{-1},$ it follows that $1/\Gamma \in L_2[-\pi,\,\pi]$, thus, $\psi^I(x) = \sum_{k\in\mathbb{Z}}\beta_k \psi(x + k),$ where $(\beta_k)_k \in l_2$ are the Fourier coefficients of $1/\Gamma.$ Therefore, $\psi^I \in W_0(\psi).$ Thus, $W_0(\psi) = W_0(\psi^I)$, so taking into account scaling, we get  $W_j(\psi) = W_j(\psi^I)$ for $j \in \mathbb{Z}$, that is, the function $\psi$ generates the same wavelet spaces $W_j(\psi)$ as the orthogonal wavelet function $\psi^I.$   

The system $\psi_{0,k},$ $k\in\mathbb{Z},$ is a Riesz basis in $W_0(\psi)$ with bounds $A,$ $B$, therefore
$\psi_{j,k},$ $k\in\mathbb{Z},$ forms a Riesz basis in $W_j(\psi),$ $j\in \mathbb{Z}$, with the same bounds. The equality  $W_j(\psi) = W_j(\psi^I)$ implies that the spaces $W_j(\psi)$, $j\in \mathbb{Z}$, form an orthogonal decomposition of $L_2(\mathbb{R}).$ Therefore, $\psi_{j,k},$ $j, k\in\mathbb{Z},$ is a Riesz basis in $L_2(\mathbb{R})$ with bounds $A,$ $B$. \hfill $\Box$    

If the sequence $\gamma_k,$ $k\in \mathbb{Z},$ is finite, then constraints on it can be formulated in slightly  another way.  
\begin{cor}
Suppose the sequence $\gamma_k,$ $k\in \mathbb{Z},$ is finite, namely, $\gamma_k=0$ for $k<N_1,$ $k>N_2,$ $N_1<N_2,$ $N_1,$ $N_2 \in \mathbb{Z},$ and the polynomial $\displaystyle \widetilde{\Gamma}(z):= z^{\max(-N_1,0)} \sum_{k=N_1}^{N_2} \gamma_k z^k$ does not have any roots on the unit circle.   
Then there exists a wavelet Riesz basis $\psi_{j,k},$ $j,k\in\mathbb{Z}$ of $L_2(\mathbb{R})$ such that 
$$
\psi(k+1/2)=\gamma_k, \ \  k\in\mathbb{Z}.
$$ 
\end{cor}

\textbf{Proof.} It follows from the  finiteness of $\gamma_k$ that $\displaystyle \Gamma(\xi)=\sum_{k=N_1}^{N_2}\gamma_k {\rm e}^{{\rm i}k \xi}$ is a trigonometric polynomial. Therefore, the inequality $|\Gamma(\xi)|^2\leq B$ is fulfilled. The inequality $|\Gamma(\xi)|^2\ge A$ holds iff the polynomial $\displaystyle \widetilde{\Gamma}(z)$ does not have any roots on the unit circle. \hfill $\Box$

If $\widetilde{\Gamma}$ has  roots on the unit circle, then one can replace these roots by close ones lying outside the unit circle. The coefficients of a polynomial continuously depend on its roots. So, the new coefficients $\gamma'_k$ are close to the old ones $\gamma_k.$ And for $\gamma'_k$ we can derive a matched wavelet.

In Theorem \ref{Mey}, we  find a Meyer wavelet that interpolates wide class of data on the lattice $1/2+3 n$, $n\in\mathbb{N}$.   

\begin{teo}
\label{Mey}
Suppose $\gamma_k \in \mathbb{R},$ $k\in\mathbb{Z}_+,$ is a sequence such that $\gamma_k = O(k^{-2-\varepsilon}),$ 
$\varepsilon>0,$   
\begin{gather}
\label{gam1}
\sqrt{2} \gamma_0 + \sqrt{2} \sum_{k=0}^{\infty} \sum_{q=0}^{\infty} \left(1+4(-1)^{q-1}2^q\right) \gamma_{2^q(2k+1)}=1,
 \\
\label{gam2}
\sqrt{2} \gamma_0 + \sqrt{2} \sum_{k=0}^{\infty} \sum_{q=0}^{\infty} \left(1-2(-1)^{q-1}2^q\right) \gamma_{2^q(2k+1)}=\frac{\sqrt{2}}{2},\\
\sqrt{2} |\gamma_0| + 3\sqrt{2} \sum_{k=0}^{\infty} \sum_{p=0}^{\infty} 
\left|\sum_{q=0}^{\infty} (-1)^{q} 2^q \gamma_{2^{q+p}(2k+1)}\right|\leq 1.
\label{gam3} 
\end{gather}
 Then there exists a Meyer wavelet $\psi^M$ such that 
$$
\psi^M\left(\frac{1}{2}+3 k\right) = \gamma_k \mbox{ for any } k \in \mathbb{N}.
$$  
\end{teo}

\textbf{Proof.} 
By definition of the Meyer wavelet, and evenness of the function $\lambda$,     we obtain
$$
\psi^M\left(\frac{1}{2}+3k\right) = \frac{1}{2\pi} \int_{\mathbb{R}} \cos(3 k \xi) \sin\lambda(\xi)\,d\xi
$$
$$
=\frac{1}{\pi} \int_{\frac{2\pi}{3}}^{\frac{4\pi}{3}} \cos(3 k \xi) \sin\left(\frac{\pi}{4}+\theta(\xi-\pi)\right)\,d\xi + 
\frac{1}{\pi} \int_{\frac{4\pi}{3}}^{\frac{8\pi}{3}} \cos(3 k \xi) \sin\left(\frac{\pi}{4}-\theta\left(\frac{\xi}{2}-\pi\right)\right)\,d\xi. 
$$
Substituting $\xi'$ for $\xi-\pi$ in the first integral and $\xi'$ for $\xi/2-\pi$ in the second one, and then using trigonometric formulas, we continue
 $$
\frac{1}{\pi} \int_{-\frac{\pi}{3}}^{\frac{\pi}{3}} (-1)^k \cos(3 k \xi') \sin\left(\frac{\pi}{4}+\theta(\xi')\right)\,d\xi' + 
\frac{1}{\pi} \int_{-\frac{\pi}{3}}^{\frac{\pi}{3}} 2\cos(6 k \xi') \sin\left(\frac{\pi}{4}-\theta\left(\xi'\right)\right)\,d\xi' 
$$
$$
=\frac{1}{\pi \sqrt{2}} \int_{-\frac{\pi}{3}}^{\frac{\pi}{3}} (-1)^k \cos(3 k \xi') \left(\sin\theta(\xi')+\cos\theta(\xi')
\right) + 2\cos(6 k \xi') \left(\cos\theta(\xi')-\sin\theta(\xi')\right)\,d\xi'
$$
$$
=\frac{1}{\pi \sqrt{2}} \int_{-\frac{\pi}{3}}^{\frac{\pi}{3}} \left((-1)^k \cos(3 k \xi') - 2\cos(6 k \xi')\right) \sin\theta(\xi')+
\left((-1)^k \cos(3 k \xi') + 2\cos(6 k \xi')\right) \cos\theta(\xi')\,d\xi'. 
$$
Since the function $\left((-1)^k \cos(3 k \xi') - 2\cos(6 k \xi')\right) \sin\theta(\xi')$ is odd it follows that 
$$
\int_{-\frac{\pi}{3}}^{\frac{\pi}{3}} \left((-1)^k \cos(3 k \xi') - 2\cos(6 k \xi')\right) \sin\theta(\xi')\,d\xi'=0.
$$
So, we get
$$
\psi^M\left(\frac{1}{2}+3k\right) =\frac{1}{\pi \sqrt{2}} \int_{-\frac{\pi}{3}}^{\frac{\pi}{3}} 
\left((-1)^k \cos(3 k \xi) + 2\cos(6 k \xi)\right) \cos\theta(\xi)\,d\xi. 
$$
Denote $h(\xi):= \cos\theta(\xi).$ The properties of the function $\theta$ implies that we need to find a function $h$ satisfying the following constraints:    
$h$ is even, $h(0)=1,$ $h(\pi/3)=1/\sqrt{2},$ $h'(0)=h'(\pi/3)=0,$ and $|h|\leq 1.$   
Let us expand the function $h$ into the Fourier series 
$\displaystyle h(\xi)=\sum_{n=0}^{\infty} \hat{h}(n) \cos(3 n \xi)$ and find coefficients $\hat{h}(n)$ to satisfy the equations $\psi^M(1/2+3k) = \gamma_k,$ $k =0,1,\dots.$ Substituting the expression for $h$ into the integral,we get
$$
 \sum_{n=0}^{\infty} \hat{h}(n) \int_{-\frac{\pi}{3}}^{\frac{\pi}{3}} 
\left((-1)^k \cos(3 k \xi) + 2\cos(6 k \xi)\right) \cos\left(3 n \xi\right)\,d\xi =
 \pi \sqrt{2} \gamma_k .
$$
Calculating the integrals we obtain the following system of equations with respect to $\hat{h}(n)$
\begin{equation} 
\label{rec}
(-1)^k \hat{h}(k)+ 2\hat{h}(2k) = 3 \sqrt2 \gamma_k, \quad k=0,1,\dots.
\end{equation} 
The system has infinitely many solutions and can be solved consequently starting from $k=0$. In this case
$\hat{h}(0)=\sqrt{2} \gamma_0.$ For $n\in \mathbb{N}$ there is a unique $p\in \mathbb{Z}_+$ and $k\in \mathbb{Z}_+,$  such that $n=2^p (2k+1)$, then solving the recurrence relation (\ref{rec})  we obtain 
$$
\hat{h}(2^p (2k+1)) = 3 \sqrt{2} \left(\frac{1}{2} \gamma_{2^{p-1}(2k+1)} - \dots + (-1)^{p-1} \frac{1}{2^p} \gamma_{2k+1} \right)+ (-1)^{p-1}\frac{1}{2^p} \hat{h}(2k+1),  
$$
and $\hat{h}(2k+1)$ can be chosen arbitrarily. To simplify checking of convergence for series let us choose
\begin{equation}
\label{odd}
\hat{h}(2k+1) = 3\sqrt{2} \sum_{q=0}^{\infty} (-1)^{q-1} 2^q \gamma_{2^q(2k+1)}, \quad  k \in \mathbb{Z}_+. 
\end{equation}
This choice is prompted by the solution of finite dimensional analogue of system (\ref{rec}).    
Then 
$$
\hat{h}(2^p(2k+1)) = 3\sqrt{2} \sum_{q=0}^{\infty} (-1)^{q} 2^q \gamma_{2^{q+p}(2k+1)}, \quad  p \in \mathbb{N}. 
$$  
Taking into account the assumption  $\displaystyle \gamma_n = O (n^{-2-\varepsilon}),$ $\varepsilon >0$, we see that the series defining  $\widehat{h}(n)$ are convergent.

Now we need to check the conditions $h(0)=1,$ $h(\pi/3)=1/\sqrt{2},$ $h'(0)=0,$ $h'(\pi/3)=0,$ and $|h(\xi)|\leq 1$. First, we prove that the series $\displaystyle  \sum_{k=0}^{\infty} k \bigl|\hat{h}(k)\bigr|$ is convergent, in this case the function $h$ is well-defined and its derivative can be calculated by term by term differentiation of the series.  
It is well known that to this end it is sufficient to prove the convergence of the series 
$\displaystyle  \sum_{k=0}^{\infty} \sum_{p=0}^{\infty}2^p (2k+1) \bigl|\hat{h}(2^p (2k+1))\bigr|.$
Using the assumption $\displaystyle \gamma_n = O (n^{-2-\varepsilon}),$ $\varepsilon >0$ we obtain 
$$
 \sum_{k=0}^{\infty} \sum_{p=0}^{\infty}2^p (2k+1) \bigl|\hat{h}(2^p (2k+1))\bigr| 
\leq 
3\sqrt{2} \sum_{k=0}^{\infty} \sum_{p=0}^{\infty}2^p (2k+1) \sum_{q=0}^{\infty} 2^q 
\bigl|\gamma_{2^{p+q}(2k+1)}\bigr| 
$$
$$
\leq C \sum_{k=0}^{\infty} \sum_{p=0}^{\infty}2^p (2k+1) \sum_{q=0}^{\infty} \frac{2^q}{2^{(p+q)(2+\varepsilon)} (2k+1)^{2+\varepsilon}}
 = C  \sum_{k=0}^{\infty} \frac{1}{(2k+1)^{1+\varepsilon}} \sum_{p=0}^{\infty}  \sum_{q=0}^{\infty}
\frac{1}{2^{(p+q)(1+\varepsilon)}}.
$$
The last series is obviously convergent. Since the series  $\displaystyle  \sum_{k=0}^{\infty} k \bigl|\hat{h}(k)\bigr|$ is convergent, it follows that 
$\displaystyle h'(\xi) =-3 \sum_{k=0}^{\infty} k \hat{h}(k)\sin 3 k \xi.$ So, the conditions $h'(0)=0,$ $h'(\pi/3)=0$ are fulfilled.

Now we turn to the conditions  $h(0)=1,$ $h(\pi/3)=1/\sqrt{2}.$
Summing (\ref{rec}) over $k=0,\dots,$ we have 
$$
\sum_{k=0}^{\infty}(-1)^k \hat{h}(k)+ 2\sum_{k=0}^{\infty}\hat{h}(2k) = 3 \sqrt2 \sum_{k=0}^{\infty}\gamma_k.
$$    
Dividing the first sum into two parts we  get
$$
-\sum_{k=0}^{\infty} \hat{h}(2k+1)+ \sum_{k=0}^{\infty}\hat{h}(2k) +\sum_{k=0}^{\infty}\hat{h}(2k) = 3 \sqrt2 \sum_{k=0}^{\infty}\gamma_k.
$$
Adding consequently $\displaystyle 4\sum_{k=0}^{\infty} \hat{h}(2k+1)$ and $\displaystyle -2\sum_{k=0}^{\infty} \hat{h}(2k+1)$ to both sides of the last equality, we obtain 
\begin{gather*}
3\sum_{k=0}^{\infty}\hat{h}(k) = 3 \sqrt2 \sum_{k=0}^{\infty}\gamma_k + 4 \sum_{k=0}^{\infty}\hat{h}(2k+1),\\
3\sum_{k=0}^{\infty} (-1)^k\hat{h}(k) = 3 \sqrt2 \sum_{k=0}^{\infty}\gamma_k -2  \sum_{k=0}^{\infty} \hat{h}(2k+1)
\end{gather*}
Substituting (\ref{odd}) and taking into account (\ref{gam1}) and (\ref{gam2}), we obtain 
 $$
h(0)=\sum_{k=0}^{\infty} \hat{h}(k) = \sqrt{2} \sum_{k=0}^{\infty}\gamma_k + 4 \sqrt{2} 
\sum_{k=0}^{\infty} \sum_{q=0}^{\infty} (-1)^{q-1} 2^q  \gamma_{2^q(2k+1)}
= 1, 
$$ 
 $$
h\left(\frac{\pi}{3}\right)=\sum_{k=0}^{\infty}(-1)^k \hat{h}(k) = \sqrt{2} \sum_{k=0}^{\infty}\gamma_k -2 \sqrt{2} 
\sum_{k=0}^{\infty} \sum_{q=0}^{\infty} (-1)^{q-1} 2^q  \gamma_{2^q(2k+1)}
= \frac{\sqrt{2}}{2}. 
$$
The last condition we need to check, $|h(\xi)|\leq 1,$ follows from (\ref{gam3})
$$
|h(\xi)|\leq \sum_{n=0}^{\infty} |\widehat{h}(n)| = 
\sqrt{2} |\gamma_0| + 3\sqrt{2} \sum_{k=0}^{\infty} \sum_{p=0}^{\infty} 
\left|\sum_{q=0}^{\infty} (-1)^{q} 2^q \gamma_{2^{q+p}(2k+1)}\right|\leq 1. \qquad \Box
$$

\section*{Acknowledgments}
The author is supported by the
Russian Science Foundation under grant No. 18-11-00055

\end{document}